%% file: Poisson_Debiased.tex
\documentclass{article}
\input{arxiv_style}
\usepackage{authblk}

\begin{document}
\title{Debiased LASSO under Poisson-Gauss Model}
\author[1]{Pedro Abdalla}
\author[2]{Gil Kur}
\affil[1]{Department of Mathematics, ETH Z\"{u}rich}
\affil[2]{Department of Computer Science, ETH Z\"{u}rich}
\maketitle

\begin{abstract}
Quantifying uncertainty in high-dimensional sparse linear regression is a fundamental task in statistics that arises in various applications. One of the most successful methods for quantifying uncertainty is the debiased LASSO, which has a solid theoretical foundation but is restricted to settings where the noise is purely additive. Motivated by real-world applications, we study the so-called Poisson inverse problem with additive Gaussian noise and propose a debiased LASSO algorithm that only requires $n \gg s\log^2p$ samples, which is optimal up to a logarithmic factor.
\end{abstract}

\section{Introduction}
Over the past two decades, regularized estimators, such as the LASSO, have been extensively employed to address high-dimensional regression problems where the number of parameters exceeds the number of observations. In this context, most research has focused on the point estimation of the regression vector. More accurately, they consider the following random sparse regression model: given $n$ i.i.d pairs $(y_1,a_1),\ldots,(y_n,a_n)$ with $y_i\in \mathbb{R}$ and $a_i\in \mathbb{R}^p$ ($n\ll p$) contaminated with an additive Gaussian noise $\sigma g_i \sim N(0,\sigma^2)$
\begin{equation}
\label{eq:basic_regression}
    y_i = \langle a_i,x^{\ast}\rangle + \sigma g_i,
\end{equation}
estimate the unknown vector $x^{\ast}\in \mathbb{R}^p$, assuming that $x^{\ast}$ is sparse. One of the most common sparsity assumptions is that $x^{\ast}$ only has $s$ non-zero entries, i.e, $x^{\ast}$ is an $s$-sparse vector.

Arguably, the most prominent method to solve \eqref{eq:basic_regression} is the LASSO that consists in the following optimization program: for some well-chosen tuning parameter $\lambda$ and a matrix $A = [a_1,\ldots,a_n]^T \in \mathbb{R}^{n\times p}$, solve
\begin{equation}
\label{eq:LASSO_optimization}
    \widehat{x}(A,\lambda,y) = \argmin_{x \in \mathbb{R}^p} \|y-Ax\|_2^2 +\lambda\|x\|_1.
\end{equation}
Given some certain technical conditions on design matrix $A$, the well established LASSO theory \cite*{wainwright2019high,buhlmann2011statistics}, asserts that the LASSO estimator enjoys (sharp) consistency
\begin{equation}
\label{eq:convergence)_rate_standard_LASSO}
    \|\widehat{x}-x^{\ast}\|_2 \lesssim \sigma\sqrt{\frac{s\log p}{n}},
\end{equation}
where $\lesssim$ (or $O(\cdot)$) denotes inequality up to a multiplicative \emph{absolute} constant $C>0$, i.e., a constant that is independent of any other parameter of the problem, including $s,n,p$.

As it often occurs, consistency alone does not suffice in many applications; one is interested in constructing confidence intervals to quantify the uncertainty of the point estimation of the regression vector. Indeed, quantifying the uncertainty of the LASSO estimator is crucial for more guided decision-making, obtaining $p$-values, and constructing hypothesis testing for variable selection \cite*{javanmard2014confidence,van2014asymptotically}. Moreover, there has been recent interest in quantifying uncertainty in medical imaging systems to ensure reliable medical diagnosis tasks in MRI \cite*{hoppe2023uncertainty,hoppe2023sampling}.

Mathematically, the goal is to construct confidence intervals $I_{i}(\alpha)$ that satisfies (as $n$ goes to infinity)
\begin{equation}
\label{eq:asymptotically_valid_CI}
    \mathbb{P}(x_i^{\ast} \in I_{i}(\alpha))\ge 1-\alpha-o(1),
\end{equation}
where $g(n)=o(f(n))$ (or $f(n) \gg g(n)$)  means that $\lim_{n\rightarrow \infty} g(n)/f(n) = 0$, for any two functions $f(n)$ and $g(n)$.

A key difficulty in quantifying the uncertainty of estimating $x^{\ast}$ is that the LASSO estimator does not enjoy a tractable asymptotic distribution as explained in \cite*{van2014asymptotically}. To bypass this issue, the debiased LASSO was proposed \cite*{van2014asymptotically,javanmard2018debiasing} as a debiasing method that constructs confidence intervals for the high-dimensional regression model \eqref{eq:basic_regression}. For the sake of discussion, we assume, for instance, that the random design matrix $A$ is isotropic, precisely $\mathbb{E}[A^TA] = I_{p\times p}$. The debiased method is described as follows: given the solution to the LASSO $\widehat{x}$, construct the debiased estimator $\widehat{x}^d$ by computing
\begin{equation}
\label{eq:basic_debiased}
    \widehat{x}^d:= \widehat{x} + \frac{1}{n}A^T(y-A\widehat{x}).
\end{equation}
A general strategy followed by \cite*{van2014asymptotically,javanmard2018debiasing} to analyze $\widehat{x}^d$ is to decompose $\sqrt{n}(\widehat{x}^d - x^{\ast})$ as
\begin{equation}
\label{eq:decomposition_deabised_basic}
\sqrt{n}(\widehat{x}^d - x^{\ast}) = \underbrace{\frac{\sigma}{\sqrt{n}}A^T g}_{\text{Gaussian term}} + \underbrace{\sqrt{n}(A^TA-I_p)(\widehat{x}-x^{\ast})}_{\text{reminder term}}.
\end{equation}
One key result from \cite*{javanmard2018debiasing} is that the reminder term is asymptotically vanishing whenever $n \gg s\log^2p$, therefore conditioned on $A$, the random variable $\sqrt{n}(\widehat{x}_i^d-x_i^{\ast})$ converges to a centered Gaussian random variable with variance $\sigma^2(A^TA)_{i,i}$. It remains to construct an honest confidence interval for this Gaussian term. To this end, let $\Phi(x)$ be the cumulative density function of a standard Gaussian and set $\delta_i:=\sigma/\sqrt{n}(A^TA_{ii})^{1/2}$ (the standard deviation of the Gaussian term). The estimated confidence interval for $x_i^{\ast}$ is given by
\begin{equation}
\label{eq:confidencen_levels}
\widehat{I}_i(\alpha) = \widehat{x}_i^d \pm \Phi^{-1}(1-\alpha/2)\delta_i.
\end{equation}
It follows immediately that $\widehat{I}_i(\alpha)$ satisfies \eqref{eq:asymptotically_valid_CI}.

On the other hand, $\sigma$ may not be known, and then we must estimate it to compute the covariance matrix of the Gaussian term. It suffices to construct an estimator $\widehat{\sigma}$ that converges to $\sigma$. Fortunately, there are methods available in the literature that estimate the standard deviation of the noise $\sigma$. One primary example is the so-called scaled LASSO \cite*{sun2012scaled}
\begin{equation}
\label{eq:scaled_LASSO}
    \{\widehat{x}(A,\lambda,y),\widehat{\sigma}(\lambda)\} = \argmin_{x \in \mathbb{R}^p,\sigma>0} \frac{1}{n\sigma}\|y-Ax\|_2^2 + \frac{\sigma}{2} +\lambda\|x\|_1,
\end{equation}
that satisfies $\widehat{\sigma}\rightarrow \sigma$ in probability for some well-chosen $\lambda$. Thus the general strategy employed in \cite*{javanmard2018debiasing} is to use $\widehat{\sigma}A^{T}g$ to quantify the uncertainty of the estimator $\widehat{x}^{d}$ around its mean $x^{\ast}$. We remark that the debiased LASSO method is optimal in the following sense. Up to a log factor, it has the best possible dependency on $n$ because $n\gg s\log p$ is necessary for the LASSO to be consistent \cite*{wainwright2019high}. Also, minimax results from \cite*{javanmard2018debiasing} establish that the length of the confidence intervals is optimal up to an absolute constant. 

However, in many real-world applications, the noise component in \eqref{eq:basic_regression} is not only additive. Indeed, Poisson noise arises in a variety of settings such as Positron Emission Tomography (PET) and medical imaging \cite*{willett2003platelets,lingenfelter2009sparsity,schmidt2009optimal,gravel2004method}, astronomy \cite*{starck2009astronomical}, and social media \cite*{xu2012socioscope}.  In these settings, the observations $y_i$ are characterized by counting the number of events modeled by a Poisson process.  For example, in $X$-ray films, the image is created by the accumulation of photons in the detector and the observations are modeled by a Poisson distribution that counts those photons (see \cite{gravel2004method} for a more accurate description of photon-limited imaging and related models with Poisson noise). In addition to the Poisson noise, an additive Gaussian noise is present in medical imaging applications due to the thermal fluctuations.

Motivated by the applications above, we study the following model (which we refer to as the Poisson-Gauss model): We design a sampling matrix $A\in \{0,1\}^{n\times p}$ and collect $n$ i.i.d pairs $(y_1,a_1),\ldots, (y_n,a_n)$ with $y_i \in \mathbb{R}$ and $a_{i}\in \{0,1\}^p$ that satisfy
\begin{equation}
\label{eq:main_regression_model}
    y_i \sim \mathcal{P}((Ax^{\ast})_i) + \sigma g_i,
\end{equation}
where $\mathcal{P}(\lambda)$ is a Poisson random variable with parameter $\lambda$ i.e, for every integer $k \geq 0$
\begin{equation*}
    \mathbb{P}(\mathcal{P}(\lambda)=k)  := \frac{\lambda^k e^{k}}{k!}.
\end{equation*}
Our goal is to estimate a \emph{non-negative} $s$-sparse vector $x^{\ast}$ based on the observations $(y_1,a_1),\ldots,(y_n,a_n)$. 
The first observation is we must consider non-negative vectors $x^{\ast} \geq 0$, otherwise the Poisson distribution does not make sense. Secondly, notice that now the noise is correlated with the design matrix $A=[a_1,\ldots,a_n]^T$, in particular, the variance of $y_i$ now depends on the value of $\langle a_i,x^{\ast}\rangle$, i.e., it presents heteroscedasticity. In contrast, in \eqref{eq:basic_regression} the variance of $y_i$ is $\sigma^2$ for all $i\in [n]$ (homocedastic).

We highlight that even the Poisson noise alone calls for a different machinery to solve \eqref{eq:main_regression_model}, even without Gaussian additive noise, when $\sigma =0$. Indeed, a line of work \cite*{hunt2018data,bohra2019variance,motamedvaziri2013sparse,raginsky2010compressed} investigated the question of how to construct consistent estimators for $x^{\ast}$ based on different assumptions on the design matrix $A$. In this paper, we study the case of $A$ generated randomly by $0/1$ Bernoulli i.i.d entries with an \emph{absolute} constant parameter $q \in (0,1)$. Namely,
$A$ has i.i.d entries $a_{ij}$, that each is independently distributed as 
\begin{equation}
\label{eq:design_matrix_A}
    a_{ij}= 
    \begin{cases}
        1, \quad \text{with probability $q$}\\
        0, \quad \text{otherwise},
    \end{cases}
\end{equation}
and for simplicity of notation, we denote by $A \sim Bern(q)$.  For example, this model appears in photon-limited imaging \cite*{hunt2018data}. Remarkably, even in the simpler case when $\sigma=0$, the analysis of the LASSO type estimator is delicate. The authors in \cite*{hunt2018data} showed that the LASSO estimator admits a suitable modification that leads to an accurate estimate of $x^{\ast}$. Specifically, by considering the normalized (and centered) version of $A$, i.e. 
\begin{equation}
\label{eq:tilde_A}
\Tilde{A}:= \frac{A}{\sqrt{nq(1-q)}} - \frac{q\mathbf{J}_{n\times p}}{\sqrt{nq(1-q)}},
\end{equation}
where $\mathbf{J}_{n\times p}$ is the $n\times p$ all-ones matrix,  and adjusting the response vector $y$ to
\begin{equation}
\label{eq:tilde_y}
\Tilde{y}:= \frac{1}{(n-1)\sqrt{nq(1-q)}}\left(ny - \mathbf{1}_{n\times 1}\sum_{i=1}^n y_i\right),
\end{equation}
where $\mathbf{1}_{n\times 1}$ is the all-ones $n \times 1$ vector. They showed that there is a well-chosen data-dependent choice of regularization parameter $\lambda'$ for which solution of the LASSO $\widehat{x}(\Tilde{A},\lambda',\Tilde{y})$ satisfies the following: in the range $n>s\log p$ and $\|x^{\ast}\|_{1}=O(n/\log p)$, we have that
\begin{equation*}
    \|\widehat{x}- x^{\ast}\|_2 = O\left( \sqrt{\frac{\|x^{\ast}\|_1s\log p}{nq}}\right)
\end{equation*}
with a high probability -- in this paper, it means that with probability tending to one as $n$ goes to infinity. Moreover, the authors showed that this rate is minimax optimal for an absolute constant $q \in (0,1)$. However, it be improved via a weighted LASSO algorithm (we refer to \cite*{hunt2018data} for more details). We do not treat this case here as we consider $q$ to be an absolute constant. The scaling with $\sqrt{\|x^{\ast}\|_1}$ comes from the fact that the variance of $y$ in \eqref{eq:main_regression_model} scales with $\|x^{\ast}\|_1$ so to compare with the standard convergence rate of the LASSO, we need to consider the normalized error $\|\widehat{x}-x^{\ast}\|_2/\sqrt{\|x^{\ast}\|_1}$.

Finally, we remark that we are not aware of any \emph{efficient} estimator, i.e., an estimator that can be computed in polynomial runtime, that achieves the same rate of convergence when $\sigma >0$, i.e., under the presence of an additive Gaussian noise as well. A preliminary step of our analysis is to extend some of the results in \cite*{hunt2018data} to the general Poisson-Gauss setup. Furthermore, to the best of our knowledge, there is no method to quantify the uncertainty under Poisson noise for high-dimensional sparse regression. In particular, it is unclear if the method \eqref{eq:basic_debiased} leads to tractable distributions under Poisson noise alone. The main goal of this manuscript is to fulfill this gap in the literature by constructing an estimator $\widehat{x}$, such that $\widehat{x}|_{A}$ (i.e. the conditional law of $\widehat{x}$ on $A$) converges entrywise (in distribution) to a Gaussian with mean $x^{\ast}$ and known variance. In what follows, $\underset{D}{\rightarrow}$ stands for convergence in distribution. Our main result is the following.

\begin{theorem}[Main Result]
\label{thm:main_result}
Consider the model of \eqref{eq:main_regression_model} above with a design matrix $A \sim Bern(q)$ for some absolute constant $q\in (0,1/2)$, and assume that there is a sufficiently large  constant $C > 0$ for which $\|x^{\ast}\|_{1}\le Cn/\log p$. Then, in the regime of $n,p\rightarrow \infty$ and $ s\log^2 p/n \rightarrow 0$, there exists an efficient estimator $\widehat{x}^d$ that satisfies
\begin{equation*}
    \sqrt{\frac{n}{\|x^{\ast}\|_1}}(\widehat{x}_i^d-x_i^{\ast})|_{A} \underset{D}{\rightarrow} N(0,\sigma(q)^2),
\end{equation*}
where $\sigma(q) = (1/\sqrt{q(1-q)})\sqrt{1+(\sigma^2\Tilde{A}^T\Tilde{A}/(n\|x^{\ast}\|_1))}$. Furthermore, when $\sigma$ or $\|x^{\ast}\|_1$ are unknown,  there are two efficient estimators $\widehat{\mu}$ and $\widehat{\nu}$ that converges to $\|x^{\ast}\|_1$ and $\sigma^2/\|x^{\ast}\|_1$ almost surely, respectively.
\end{theorem}

We remark that the normalization factor $1/\sqrt{\|x^{\ast}\|_1}$ is only needed if $\|x^{\ast}\|_1$ diverges with $n$ or $p$. It is usually assumed to be a constant and known, but we do not need it here. Indeed, we only require that $\|x^{\ast}\|_1=O(n/\log p)$ which is mild and it holds in the practically relevant case when $\|x^{\ast}\|_{\infty}=O(1)$ because $x^{\ast}$ is $s$-sparse and $n\gg s\log p$ is necessary.

Non-asymptotic probability estimates and convergence rates follow from the proof of the main result. The optimality of our result is similar to the one from \cite*{javanmard2018debiasing}; the dependence on $n$ is optimal up to a log factor, and the length of the interval is optimal up to an absolute constant. Finally, we believe that using our techniques can extend our results to the case when $x^{\ast}$ is sparse on a dictionary basis. 

The rest of this paper is organized as follows. In Section \ref{sec:UQ_poisson_alone}, we show how to quantify uncertainty under the Poisson Noise alone. In Section \ref{sec:consistency_poissongauss}, we extend some results about the consistency of the LASSO for the Poisson-Gauss setup. In Section \ref{sec:main_result}, we wrap up all the theoretical results into a final estimator. Section \ref{sec:numerical} is dedicated to numerical experiments. 

Before we move to the next section, we introduce some notation: $X|_{Y}$ refers to the random variable $X$ conditioned on $Y$,  and $\underset{P}{\to}$ means convergence in probability. We use $C,c > 0$ for absolute constants whose values may change from line to line. When the constant depends on the parameter $q$, we write $C(q)$.  For $p\in [1,\infty]$, the standard $\ell_p$ norm is denoted by $\|\cdot\|_p$. Finally, we drop the indices in $\mathbf{J}_{n\times p},\mathbf{1}_{n\times 1}$ when it is clear from the context.

\section{Uncertainty Quantification under Poisson Noise}
\label{sec:UQ_poisson_alone}
In this section, we tackle the problem of quantifying the uncertainty of the LASSO estimator under Poisson noise in the simplified case where there is no Gaussian additive noise. Namely, we prove Theorem \ref{thm:main_result}, when $\sigma=0$ in \eqref{eq:main_regression_model}.

To start, let us recall the model in the case of interest here
\begin{equation}
\label{eq:main_regression_model_poisson_only}
    y_i \sim \mathcal{P}((Ax^{\ast})_i),
\end{equation}
and recall that $A \sim Bern(q)$ for some fixed constant $q \in (0,1)$. As mentioned briefly in the introduction, the starting point of our analysis is the following result from \citep[Proposition 4]{hunt2018data}. To describe it accurately, we first require some notation. For the rest of this work, we assume that $q \in [0,1/2]$. Next, we define for every $k\in [p]$, the vector $V_k \in \mathbb{R}^n$ given entrywise by
\begin{equation*}
V_{k,l}:= \left(\frac{na_{l,k}-\sum_{i=1}^n a_{ik}}{n(n-1)q(1-q)}\right)^2.
\end{equation*}
To reduce the notation burden, define 
\begin{equation*}
W:= \max_{u,k\in [p]} \sum_{i=1}^n a_{i,u} V_{k,i}.
\end{equation*}
Next, we define an estimator for $\|x^{\ast}\|_1$ via
\begin{equation}
\label{eq:estimator_Nhat}
\widehat{N}:= \frac{\left(\sqrt{1.5\log p} + \sqrt{2.5 \log p + \sum_{i=1}^ny_i}\right)^2}{\sqrt{nq-\sqrt{6nq(1-q)\log p}-(1-q)\log p}},
\end{equation}
where the explicit constants in $\widehat{N}$, and later in the estimator $\widehat{d}$ below, arrive from the analysis of the LASSO estimator and do not have any special meaning. Next, we define 
$\hat{d}$, via
\begin{equation}
\label{eq:weight_d_LASSO}
\widehat{d} = \sqrt{6\widehat{N}W\log p} + \frac{\log p}{(n-1)q(1-q)} + \frac{378  \log p}{n}\left(1+ \frac{1-q}{q}\left(\frac{3\log p}{n}\right)\right)\widehat{N}.
\end{equation}
 We are ready to state their result. 
\begin{theorem}
\label{thm:consistency_LASSO_Poisson_only}
Consider the sparse regression model \eqref{eq:main_regression_model} with $\sigma=0$. Assume that $\|x^{\ast}\|_{1}\le C n/\log p$ for some absolute constant $C > 0$. Let $\Tilde{A}$, $\Tilde{y}$ and $\widehat{d}$ is given by \eqref{eq:tilde_A}, \eqref{eq:tilde_y} and \eqref{eq:weight_d_LASSO}, respectively. Then for any $\gamma>2$, the solution $\widehat{x}(\Tilde{A},\gamma \widehat{d},\Tilde{y})$ of \eqref{eq:LASSO_optimization} satisfies, with probability at least $1-e^{-n}-p^{-1}$, that
\begin{equation*}
    \|\widehat{x}-x^{\ast}\|_2 \le \kappa(\gamma) \sqrt{\frac{\|x^{\ast}\|_1s\log p}{nq}},
\end{equation*}
where  $\kappa(\gamma)$is a constant that depends only on $\gamma$.
\end{theorem}
Intuitively, the choice of $\Tilde{A}$ and $\Tilde{y}$ is to make the coherence type parameter $\|\Tilde{A}^T(\Tilde{y}-\Tilde{A}x^{\ast})\|_{\infty}$ as small as possible, in particular, smaller than the regularization parameter $\gamma\widehat{d}$. It follows from the theory of LASSO that when this is the case, the LASSO estimator has some guarantees. We will see a precise statement of this type in the next chapter. 

Next, as discussed in \cite*{hunt2018data}, the convoluted expression for choice of $\widehat{d}$ in Theorem \ref{thm:consistency_LASSO_Poisson_only} could be avoided under the assumption that $\|x^{\ast}\|_1$ is known. This is due to the fact that one would like to replace $\widehat{d}$ by 
\begin{equation*}
    d=\sqrt{\frac{\|x^{\ast}\|_1 \log p }{nq(1-q)}}.
\end{equation*}

In the next section, we will present an alternative way to construct an estimator $\widehat{\mu}$ for $\|x^{\ast}\|_1$. Equipped with such an estimator, one could set 
\begin{equation}
\label{eq:easy_choice_regularization}
\widehat{d} = \sqrt{\frac{\widehat{\mu} \log p }{nq(1-q)}}
\end{equation}
as a choice of regularization parameter for the LASSO.

Next, equipped with Theorem \ref{thm:consistency_LASSO_Poisson_only}, we are now ready to start the construction of our proposed debiasing method. From the theoretical viewpoint, it will be necessary to create independence between our design matrix $A$ and the LASSO estimator $\widehat{x}$. More accurately, without loss of generality, assume that we collect $2n$ samples. We split the sample into two batches: the first one is $(y_1,a_1),\ldots,(y_n,a_n)$ which we denote by $(y,A)$ and the second one $(y_{n+1},a_{n+1}),\ldots,(y_{2n},a_{2n})$ which we denote by $(y',A')$. From the practical perspective, data splitting is not interesting as it discards half of the sample, and in many experiments, it seems unnecessary. As discussed in \cite{javanmard2018debiasing}, it is an open problem to analyze the performance of the debiased LASSO for non-Gaussian design matrices even in simplest the case \eqref{eq:basic_regression} of an absence of Poisson noise. In the same paper, the authors used a delicate leave-one-out argument that requires exclusive properties of Gaussian design matrices.

In what follows, we first focus on the case when $\|x^{\ast}\|_1$ is known to facilitate the analysis. We relax this assumption later. For instance, set
\begin{equation}
\label{eq:expression_known_bias}
B :=(\|x^{\ast}\|_1-\|\widehat{x}\|_1)\frac{q}{\sqrt{nq(1-q)}}\widetilde{A}^T \mathbf{1},
\end{equation}
our algorithm is described below.
\begin{algorithm}
\caption{Debiased Algorithm Under Poisson with Known $\|x^{\ast}\|_{1}$}
\label{Alg:debiased_estimator}
\begin{algorithmic}
\Require Design matrix $ \begin{pmatrix}
A\\
A'
\end{pmatrix} \in \{0,1\}^{2n\times p}$, a response vector $(y,y')\in \mathbb{R}^{2n}$ from \eqref{eq:main_regression_model_poisson_only}. 
\State
\State Compute $\Tilde{A}$, $\Tilde{A'}$ by \eqref{eq:tilde_A} and $\Tilde{y}'$ by \eqref{eq:tilde_y}.

\State Compute the LASSO solution $\widehat{x}$ as in Theorem \ref{thm:consistency_LASSO_Poisson_only} using $\Tilde{A'}$ and $\Tilde{y}'$.

\State Compute the known bias term $B$ using \eqref{eq:expression_known_bias}. 
\\

\Return
\begin{equation}
\label{eq:debiased_estimator}
    \widehat{x}^{d}:= \widehat{x} + \frac{1}{\sqrt{nq(1-q)}}\Tilde{A}^T(y-A\widehat{x}) - \frac{1}{\sqrt{n}}B.
\end{equation}
\end{algorithmic}
\end{algorithm}

Recall that $\mathbf{J}=\{1\}^{n\times p}$, the starting point of our analysis is the following decomposition
\begin{equation}
\label{ineq:steps_debiasing_poisson_only}
\begin{split}
&\sqrt{n}({\widehat{x}}^{d}-x^{\ast})\\
&= \sqrt{n}(\widehat{x}-x^{\ast}) + \frac{1}{\sqrt{q(1-q)}}\widetilde{A}^T(\mathcal{P}(Ax^{\ast})-Ax^{\ast}) + \frac{1}{\sqrt{q(1-q)}}\widetilde{A}^T A(x^{\ast}- \widehat{x}) - B\\
&=\sqrt{n}(\widehat{x}-x^{\ast}) + \frac{1}{\sqrt{q(1-q)}}\widetilde{A}^T(\mathcal{P}(Ax^{\ast})-Ax^{\ast}) + \sqrt{n}\widetilde{A}^T \frac{A-q\mathbf{J}}{\sqrt{nq(1-q)}}(x^{\ast}- \widehat{x}) \\
&+ \widetilde{A}^T\frac{q\mathbf{J}}{\sqrt{q(1-q)}}(x^{\ast}-\widehat{x}) - B\\
&=\frac{1}{\sqrt{q(1-q)}}\widetilde{A}^T(\mathcal{P}(Ax^{\ast})-Ax^{\ast}) + \sqrt{n}(\widetilde{A}^T \widetilde{A}-I)(x^{\ast}- \widehat{x}) + \underbrace{(\|x^{\ast}\|_1-\|\widehat{x}\|_1)\frac{q}{\sqrt{nq(1-q)}}\widetilde{A}^T \mathbf{1}}_{=B} - B\\
&= \underbrace{\frac{1}{\sqrt{q(1-q)}}\widetilde{A}^T(\mathcal{P}(Ax^{\ast})-Ax^{\ast})}_{:=\eta}+ \underbrace{\sqrt{n}(\widetilde{A}^T \widetilde{A}-I)(x^{\ast}- \widehat{x})}_{:=\Delta}.
\end{split}
\end{equation}
In analogy to the decomposition \eqref{eq:decomposition_deabised_basic}, the random vector $\eta$ is the noise term that provides the uncertainty replacing the Gaussian term, $\Delta$ is the reminder negligible term and $B$ is an additional bias term (which we remove) that comes from the fact that we used $A$ instead of $\Tilde{A}$ because $y$ is non-linear on $A$.

\subsection{The case when $\|x^{\ast}\|_1$ is known}
Our goal is to build confidence intervals for the LASSO solution $\widehat{x}(\Tilde{A'},\gamma d,\Tilde{y'})$ when the one-dimensional parameter $\|x^{\ast}\|_1$ is known. In the next section, we relax this assumption by proposing an estimator for $\|x^{\ast}\|_1$ that is arguably simpler than the construction of $\widehat{N}$ in \eqref{eq:estimator_Nhat} above.

The core of the argument is to show the following:
\begin{proposition}
\label{prop:main_chapter_2}
Under the notation above, the following holds:
\begin{enumerate}
\label{enum:key_items}
    \item The exact value of $B$ is known; in particular, it does not depend on $x^{\ast}$.
    \item Each entry of the vector $\eta/\sqrt{\|x^{\ast}\|_1}$ conditioned on $A$ converges in distribution to $N(0,1/q(1-q))$ .
    
    \item $\Delta/\sqrt{\|x^{\ast}\|_1}$ is negligible, precisely $\|\Delta/\sqrt{\|x^{\ast}\|_1}\|_{\infty} \underset{P}{\rightarrow} 0$.
\end{enumerate}
\end{proposition}
It follows that all the terms above are either known (after $A$ is drawn) or negligible. Moreover, notice that Proposition \ref{prop:main_chapter_2} implies Theorem \ref{thm:main_result} in the simplified case when $\sigma=0$. The rest of this section is dedicated to the proof of Proposition \ref{prop:main_chapter_2}. As the first item is straightforward because, as $\|x^{\ast}\|_1$ is known, it remains to prove Items 2 and 3.
\medskip

\noindent\textbf{The noise term $\eta$:}
Consider the random variable
\begin{equation*}
Z:=\frac{A^T-qJ^T}{\sqrt{nq(1-q)}}(\mathcal{P}(Ax^{\ast})-Ax^{\ast}).
\end{equation*}
We claim that each entry of $Z/\|x^{\ast}\|_1$ conditioned on $A$ converges to a centered Gaussian. Clearly, if this claim is true, then it is also true for $Z/(\|x^{\ast}\|_1\sqrt{q(1-q)})$.

The key standard fact that will enable us to establish asymptotic normality is the following result about Poisson distributions \cite*[Exercise 2.3.8]{vershynin2018high}. 
\begin{lemma}
\label{lemm:Poisson_converges_Gaussian}
There exists an absolute constant $C>0$ such that for any $t\in \mathbb{R}$ and positive integer $\lambda>0$ the following holds:
\begin{equation*}
    \left|\mathbb{P}\left(\frac{\mathcal{P}(\lambda)-\lambda}{\sqrt{\lambda}} \ge t\right)- \mathbb{P}(g\ge t)\right|\le \frac{C}{\sqrt{\lambda}},
\end{equation*}
where $g \sim N(0,1)$. In particular, it holds that
\begin{equation*}
    \lim_{\lambda\rightarrow \infty}\left|\mathbb{P}\left(\frac{\mathcal{P}(\lambda)-\lambda}{\sqrt{\lambda}} \ge t\right)- \mathbb{P}(g\ge t)\right|=0.
\end{equation*}
\end{lemma}
The proof is a direct consequence of the central limit theorem as a Poisson random variable with parameter $\lambda$ can be expressed as the sum of $\lambda$ (take a sub-sequence of integer $\lambda$'s) Poisson random variables with mean and variance one. Consequently, the convergence rate follows immediately from the Berry-Essen theorem.

We now apply Lemma \ref{lemm:Poisson_converges_Gaussian} to prove the second item in Proposition \ref{prop:main_chapter_2}.
\begin{proof}

Since the entries of $Z$ are identically distributed, it is enough to show that the first entry of $Z$, namely
\begin{equation*}
    Z_1:= \frac{1}{\sqrt{n}} \sum_{i=1}^n\frac{1}{\sqrt{q(1-q)}}(a_{1i}-q)(\mathcal{P}(\langle a_{i},x^{\ast}\rangle) - \langle a_i,x^{\ast}\rangle)
\end{equation*}
converges to a centered Gaussian distribution. To this end, let $I_{+}$ the indexes for which $a_{1i} = 1$ and $I_{-}$ the indexes for which $a_{i1}=0$. Clearly, 
\begin{equation*}
Z_{1} = Z_{1,+} + Z_{1,-},
\end{equation*}
where 
\begin{equation*}
Z_{1,+}:=\frac{1-q}{\sqrt{nq(1-q)}}\sum_{i\in I_{+}}(\mathcal{P}(\langle a_{i},x^{\ast}\rangle) - \langle a_i,x^{\ast}\rangle),
\end{equation*}
and 
\begin{equation*}
Z_{1,-}:=\frac{-q}{\sqrt{nq(1-q)}}\sum_{i\in I_{-}}(\mathcal{P}(\langle a_{i},x^{\ast}\rangle) - \langle a_i,x^{\ast}\rangle).
\end{equation*}
We first analyze $Z_{1,+}$. From the fact that the sum of Poisson random variables is also distributed as Poisson, we obtain that
\begin{equation*}
Z_{1,+} = \frac{1-q}{\sqrt{nq(1-q)}} \mathcal{P}\left(\left\langle\sum_{i\in I_{+}}a_{i},x^{\ast}\right \rangle\right) -  \left(\left\langle\sum_{i\in I_{+}}a_{i},x^{\ast}\right \rangle\right).
\end{equation*}
To avoid confusion with the regularization parameter $\lambda$, we set $\lambda_{\mathcal{P}}$ to be the Poisson parameter 
\begin{equation*}
\lambda_{\mathcal{P}}:= \left(\left\langle\sum_{i\in I_{+}}a_{i},x^{\ast}\right \rangle\right).
\end{equation*}
Multiplying and dividing $Z_{1,+}$ by $\sqrt{\lambda_{\mathcal{P}}}$, we obtain that
\begin{equation*}
Z_{1,+} = (1-q)\sqrt{\frac{\lambda_{\mathcal{P}}}{q(1-q)n}} \frac{\mathcal{P}(\lambda_{\mathcal{P}}) - \lambda_{\mathcal{P}}}{\sqrt{\lambda_{\mathcal{P}}}}.
\end{equation*}
 For notation simplicity, we define $W_i$'s as 
\begin{equation*}
\lambda_{\mathcal{P}} = \sum_{i\in I_{+}}\sum_{j=1}^pa_{ij}x^{\ast}_j := \sum_{i\in I_{+}} W_i,
\end{equation*}
and consider the function $$f(W_1,\ldots,W_{|I_{+}|}):=\sum_{i\in I_{+}}W_i.$$ Next, for every pair $(w_i,w'_i)$, we have that
\begin{equation*}
|f(W_1,\ldots, w_i,\ldots,W_{|I_{+}|}) - f(W_1,\ldots, w'_i,\ldots,W_{|I_{+}|})| \le |\sum_{j=1}^p (a_{ij}-a'_{ij})x^{\ast}_{j}| \le 2\|x^{\ast}\|_1,
\end{equation*}
thus by the bounded difference inequality \cite*[Theorem 2.9.1]{vershynin2018high}, 
\begin{equation*}
\mathbb{P}\left(|I_{+}|q\|x^{\ast}\|_1 - f(W_1,\ldots,W_{|I_{+}|})\ge t\right)\le e^{-t^2/(|I_{+}|\|x^{\ast}\|_1^2)}.
\end{equation*}
Setting $t= q|I_{+}|\|x\|_1/2$, we obtain that with probability at least $1-e^{-|I_{+}|q^2/4}$
\begin{equation*}
    \lambda_{\mathcal{P}} \ge \frac{q}{2}|I_{+}|\|x^{\ast}\|_1.
\end{equation*}
Thus, on the event that $|I_{+}| \rightarrow \infty$ as $n$ goes to infinity, we guarantee that $\lambda_{\mathcal{P}} \rightarrow \infty$ and Lemma \ref{lemm:Poisson_converges_Gaussian} holds. Since $|I_{+}|$ converges almost surely to its mean $nq$, this is indeed the case. To see this, notice that by Chernoff's small deviation inequality \cite*[Exercise 2.3.5]{vershynin2018high} inequality $$|I_{+}| =\left(1\pm \frac{1}{\sqrt{\log n}}\right)nq $$ with probability at least $1-e^{-cnq/\log n}$. Therefore, Lemma \ref{lemm:Poisson_converges_Gaussian} implies that $Z_{1,+}$ converges in distribution to a Gaussian provided that $\sqrt{\lambda_{\mathcal{P}}/n}$ converges to a finite non-zero number. This is indeed the case, first recall that
\begin{equation}
\label{eq:lambda/n}
    \sqrt{\frac{\lambda_{\mathcal{P}}}{q(1-q)n}} = \left(\frac{1}{n}\sum_{i\in I_{+}}\sum_{j=1}^p \frac{a_{ij}-q}{\sqrt{q(1-q)}}x^{\ast}_{j} + \frac{|I_{+}|}{nq(1-q)}\|x^{\ast}\|_1\right)^{1/2}.
\end{equation}
By the bounded differences inequality again, there is a constant $c(q)$ such that with probability at least $1-e^{-n^2t^2/(c(q)|I_{+}|\|x^{\ast}\|_1^2)}$,
\begin{equation*}
\frac{1}{n}\sum_{i\in I_{+}}\sum_{j=1}^p \frac{a_{ij}-q}{\sqrt{q(1-q)}}x^{\ast}_{j} \le t.
\end{equation*}
We choose $t=2\sqrt{c(q)|I_{+}|\|x^{\ast}\|_1/(n \log p)}$  to ensure that
\begin{equation*}
    \frac{1}{\sqrt{\|x^{\ast}\|_1}}\sqrt{\frac{1}{n}\sum_{i\in I_{+}}\sum_{j=1}^p \frac{a_{ij}-q}{\sqrt{q(1-q)}}x^{\ast}_{j}} \lesssim \frac{1}{\sqrt{\log p}} = o(1),
\end{equation*}
with probability at least $$1-e^{-4n/\log p} \ge 1-e^{-4\log p} = 1-p^{-4},$$ that converges to one. To handle the second term in \eqref{eq:lambda/n}, recall that we normalize $Z$ by $\sqrt{\|x^{\ast}\|_1}$ and notice that almost surely
\begin{equation*}
\sqrt{\frac{|I_{+}|}{nq(1-q)}} \rightarrow \sqrt{\frac{1}{1-q}},
\end{equation*}
then it also holds that almost surely,
\begin{equation*}
    (1-q)\sqrt{\frac{\lambda_{\mathcal{P}}}{nq(1-q)\|x^{\ast}\|_1}}\rightarrow \sqrt{1-q}.
\end{equation*}
We conclude that conditionally on $A$, we have the following convergence in distribution (as $n$ diverges)
\begin{equation*}
\frac{1}{\sqrt{\|x^{\ast}\|_1q(1-q)}}Z_{1,+}\underset{D}{\rightarrow}  \frac{1}{\sqrt{q}} N(0,1).
\end{equation*}
An analogous argument shows that
\begin{equation*}
\frac{1}{\sqrt{\|x^{\ast}\|_1q(1-q)}}Z_{1,-}\underset{D}{\rightarrow}  \frac{1}{\sqrt{1-q}} N(0,1),
\end{equation*}
and by independence, $Z_{1}$ converges to a centered Gaussian distribution whose variance equals to 
\[ \frac1q +\frac 1{1-q} = \frac 1{q(1-q)}.\]
\end{proof}
\medskip

\textbf{Negligible Bias $\Delta$:}
To show that $\Delta/\sqrt{\|x^{\ast}\|_1}$ is negligible, we first estimate $\|\Delta\|_{\infty}$. To this end, we need the sub-exponential version of Bernstein's inequality that we state here for the reader's convenience. Recall that the $\psi_2$ and $\psi_1$ norms of a random variable $X$ are defined as
\begin{equation*}
  \|X\|_{\psi_2}:=\inf\{t>0: \mathbb{E}e^{X^2/t}\le 2\}\quad   \|X\|_{\psi_1}:=\inf\{t>0: \mathbb{E}e^{|X|/t}\le 2\},
\end{equation*}
and $X$ is sub-Gaussian or sub-exponential if and only if $\|X\|_{\psi_2}$ or $\|X\|_{\psi_1}$ is finite, respectively. The following inequality is the sub-exponential version of the classical Bernstein's inequality \cite*[Theorem 2.8.1]{vershynin2018high}.
\begin{proposition}
\label{ineq:Bernstein_psi_1}
Let $X_1,\ldots, X_n$ be independent, mean zero, sub-exponential random variables satisfying for some $L>0$
\begin{equation*}
    \max_{i\le n}\|X_i\|_{\psi_1}\le L.
\end{equation*}
Then there is a constant $c>0$ such that for every $t>0$,
\begin{equation*}
\mathbb{P}\left(\left|\sum_{i=1}^n X_i\right|\ge t\right)\le 2 e^{{{-c\min\left\{\frac{t^2}{nL}, \frac{t}{L}\right\}}}}.
\end{equation*}
\end{proposition}
We are ready to proceed to the proof of Item 3 in Proposition \ref{prop:main_chapter_2}.
\begin{proof}
The bias $\Delta$ may also be expressed as follows
\begin{equation*}
\Delta = \sqrt{n}(\widetilde{A}^T \widetilde{A} - I)(\widehat{x}-x^{\ast}) = \frac{1}{\sqrt{n}}\sum_{i=1}^n \frac{1}{q(1-q)}(a_{i}'-q)(a_{i}'^T-q)(\widehat{x}_i-x^{\ast}_i)
\end{equation*}
By independence between $\widehat{x}$ and $A$, setting $u$ given entrywise by $$u_{i}:=\widehat{x}_i-x^{\ast},$$ we obtain that $u$ is fixed once we condition on $\widehat{x}$, and the law of $A$ remains unchanged because it is independent from both $A'$ and $y'$. Next, we show that
\begin{equation*}
\Delta_j:=\frac{1}{\sqrt{n}}\sum_{i=1}^n \frac{1}{q(1-q)}e_{j}^T(a_{i}-q)(a_{i}^T-q)u-u_j :=\frac{1}{\sqrt{n}}\sum_{i=1}^n W_i(j),
\end{equation*}
vanishes with high probability. To this end, we apply Proposition \ref{ineq:Bernstein_psi_1}: notice that conditioned on $u$, we apply \cite*[Lemma 2.7.7]{vershynin2018high} to obtain that there is a constant $C(q)$ for which
\begin{equation*}
\begin{split}
&\left\|\frac{1}{q(1-q)}e_{j}^T(a_{i}-q)(a_{i}^T-q)u\right\|_{\psi_1}\\
&\le 2 \left\|\frac{1}{\sqrt{q(1-q)}}e_{j}^T(a_{i}-q)\right\|_{\psi_2}\left\|\frac{1}{\sqrt{q(1-q)}}(a_i^T-q)u\right\|_{\psi_2}\\
&\le 2C(q) \|u\|_{2}.
\end{split}
\end{equation*}
By Proposition \ref{ineq:Bernstein_psi_1}, there is a constant $C_1(q)$ such that
\begin{equation*}
\mathbb{P}\left(\frac{1}{\sqrt{n}}\sum_{i=1}^n W_i(j) \ge t\|u\|_2 | \widehat{x}\right) \leq 2 e^{-C_1(q)t^2}.
\end{equation*}
Therefore, setting $t= 2\sqrt{\log p}/C_1(q)$ we obtain that 
\begin{equation*}
\mathbb{P}\left(\forall j\in [p]: \ \frac{1}{\sqrt{n}}\sum_{i=1}^n W_i(j) \ge \frac{2}{C_1(q)}\sqrt{\log p}\|u\|_2 | \widehat{x}\right) \le \frac{1}{p^4}.
\end{equation*}
By union bound over all $j\in [p]$, we obtain that with probability at least $1-p^{-3}$
\begin{equation*}
    \|\Delta\|_{\infty} \le C_1(q) \|\widehat{x}-x^{\ast}\|_2\sqrt{\log p}.
\end{equation*}
On that event, by Theorem \ref{thm:consistency_LASSO_Poisson_only}
\begin{equation*}
    \frac{\|\Delta\|_{\infty}}{\sqrt{\|x^{\ast}\|_1}}  \le C_1(q) \kappa(\gamma)\log p\sqrt{\frac{s}{nq}},
\end{equation*}
that vanishes provided that $n\gg s\log^2p$.
\end{proof}
The only place where we need the sample splitting is to argue that conditioned on $\widehat{x}-x^{\ast}$, the law of the design matrix does not change. It seems that in practice (see the experiments in Chapter \ref{sec:numerical}), avoiding the sample splitting is harmless.

\subsection{The case when $\|x^{\ast}\|_1 $ is unknown}
Notice that in the last section, we only required the knowledge of $\|x^{\ast}\|_1$ to remove the known bias $B$ and to normalize $\widehat{x}^{d}-x^{\ast}$. To bypass this obstacle in the case when $\|x^{\ast}\|_1$ is unknown, all we need is to find a consistent estimator for the one-dimensional quantity $\|x^{\ast}\|_1$.

To this end, we invoke some machinery from robust statistics by slightly modifying the construction of design matrix $A$, in particular, by replacing an $\theta$-fraction of the sample. More accurately, let $(a_1,\ldots a_{(1-\theta)n})$ be an independent random vector with i.i.d entries distributed as a Bernoulli $q$ and assume without loss of generality that $(1-\theta)n$ is an integer. Recall that $\mathbf{1}_{n\times 1}$ is the all-ones vector of length $n$. The new design matrix $A_{new} \in \{0,1\}^{n\times p}$ is given by
\begin{equation}
\label{eq:new_design_matrix}
A_{new}:=
\begin{pmatrix}
-a_1-\\
\vdots\\
-a_{(1-\theta)n}-\\
-\mathbf{1}_{n\times 1}-\\
\vdots\\
-\mathbf{1}_{n\times 1}-
\end{pmatrix}.
\end{equation}
The main idea is to save an $\theta n$ rows of the sample to estimate $\|x^{\ast}\|_1$ and by doing that, we obtain that for every $j\in \{(1-\theta)n+1,\ldots, n\}$,
\begin{equation*}
    y_j = \mathcal{P}_j(\|x^{\ast}\|_1).
\end{equation*}
Notice that the random variables $y_j$'s are i.i.d copies of a random variable whose mean and variance are equal to $\|x^{\ast}\|_1$. Thus, it suffices to estimate the mean of a random variable $\mathcal{P}(\|x^{\ast}\|_1)$. Thanks to the sub-exponential tails of the Poisson distribution, the empirical mean concentrates well, meaning that it is a consistent estimator for the mean. Moreover, if we add the Gaussian noise, the observations become 
\begin{equation}
\label{eq:observations_i.i.d_poisson_gauss}
    y_j = \mathcal{P}_j(\|x^{\ast}\|_1) + \sigma g_j,
\end{equation}
which also have mean $\|x^{\ast}\|_1$. It follows that the proof for the general case \eqref{eq:main_regression_model} is exactly the same. For this reason, we state our main result for the more general model \eqref{eq:main_regression_model}. As it happens, we provide a slightly sharper result by invoking standard results from the theory of median of means estimators.
\begin{proposition}
\label{prop:estimation_l1norm_x*}
Consider the regression model \eqref{eq:main_regression_model} with the design matrix $A_{new}$. Then there is a computable estimator $\widehat{\mu}$ satisfying with probability at least $1-n^{-2}$,
\begin{equation*}
|\widehat{\mu}-\|x^{\ast}\|_1| \le 8(\sqrt{\|x^{\ast}\|_1}+\sigma)\sqrt{\frac{\log n}{\theta n}}.
\end{equation*}
\end{proposition}
The proof and the description of the algorithm are postponed to the Appendix.
The main outcome is the following : any choice of $\theta \in (0,1)$ satisfying $\theta n \gg \|x^{\ast}\|_1\log n$ leads to 
\begin{equation*}
   \max\left\{\frac{\|x^{\ast}\|_1}{\widehat{\mu}},\frac{\widehat{\mu}}{\|x^{\ast}\|_1}\right\} \rightarrow 1,
\end{equation*}
almost surely. Therefore we can still use the expression for $B$ \eqref{eq:expression_known_bias} by replacing $\|x^{\ast}\|_1$ by $\widehat{\mu}$ and the same result holds thanks to Proposition \ref{prop:estimation_l1norm_x*}. It follows immediately that $\sqrt{n/\widehat{\mu}}(\widehat{x}_i^d-x_i^{\ast})$ converges in distribution as well.

\section{Consistency in the Poisson-Gauss Model}
\label{sec:consistency_poissongauss}
In the previous section, we derived a debiasing method to quantify the uncertainty of the high-dimensional regression model \eqref{eq:main_regression_model_poisson_only} with pure Poisson noise. However, as discussed in the introduction, it is also the case that the response vector $y$ is contaminated with an independent additive Gaussian noise.

Therefore, we would like to extend our method to this Poisson-Gauss model. The first step is to obtain a result analogous to Theorem \ref{thm:consistency_LASSO_Poisson_only} to the more general case \eqref{eq:main_regression_model}. To this end, it is necessary to argue that the solution of the LASSO $\widehat{x}$ still converges to the truth $x^{\ast}$ at a potentially slower convergence rate.

Notice that if $y$ is now given by \eqref{eq:main_regression_model} then the debiased estimator $\widehat{x}^d$ \eqref{eq:debiased_estimator} satisfies the same decomposition \eqref{ineq:steps_debiasing_poisson_only} except by an additional term that contains the noisy vector $\sigma g$. More accurately, the following holds
\begin{equation}
\label{eq:decomposition_deabised_final}
\sqrt{\frac{n}{\|x^{\ast}\|_1}}(x^{d}-x^{\ast}) = \frac{1}{\sqrt{\|x^{\ast}\|_1}}\left(\eta+ \Delta + \underbrace{\sigma\frac{1}{\sqrt{q(1-q)}}\Tilde{A}^T g}_{:=\eta'}\right).
\end{equation}

Thanks to Proposition \ref{prop:estimation_l1norm_x*}, the bias $B$ can be handled in the same way as before. Also, $\eta$ admits the same analysis. It remains to analyze $\Delta$ and $\eta'$. As we will see, thanks to the new consistency results (first part of Proposition \ref{main_prop_consistency} below), we can still guarantee that $\Delta/\sqrt{\|x^{\ast}\|_1}$ is negligible. Moreover, the term $\eta'$ is also Gaussian conditionally on $A$ and is independent from $\eta$. Therefore, the law of $\eta+\eta'$ becomes a centered multivariate Gaussian once we condition to $A$.

As described in the introduction, our main motivation to consider the additive Gaussian noise in the Poisson model was due to thermal fluctuations in medical imaging problems. In those cases, the variance $\sigma^2$ of the additive Gaussian noise can be assumed to be known because it can be directly measured by the system; see \cite*{gravel2004method} for more details. 

Alternatively, one could construct an estimator for $\sigma$ that suffices for our purposes. This is the main content of the second part of our main result of this section, namely Proposition \ref{main_prop_consistency} below. To state it, recall that in the design of $A_{new}$ \ref{eq:new_design_matrix} we saved $\theta n$ samples to construct an estimator for $\|x^{\ast}\|_1$.

\begin{proposition}
\label{main_prop_consistency}
Let $y,\Tilde{y}\in \mathbb{R}^n$ as in \eqref{eq:main_regression_model} and \eqref{eq:tilde_y}, respectively. Assume that $\|x^{\ast}\|_1\le Cn/\log p$ for some constant $C>0$. Set $A,\Tilde{A}$ as in \eqref{eq:design_matrix_A} and \eqref{eq:tilde_A}. For any $\gamma>2$, consider the solution to \eqref{eq:LASSO_optimization} $\widehat{x}(\Tilde{A},\lambda,\Tilde{y})$ with the choice of
\begin{equation*}
    \lambda= \gamma\left(d_{old}+ \sigma\sqrt{\frac{8\log p}{n}} + \sigma \sqrt{\frac{c_q\log p}{q(1-q)n}}\right),
\end{equation*}
where $d_{old}$ is given by \eqref{eq:easy_choice_regularization} or \eqref{eq:weight_d_LASSO}, and $c_q=(1+\sqrt{1/q(1-q)})$. Then the vector $\widehat{x}$ satisfies with high-probability
\begin{equation*}
\|\widehat{x} - x^{\ast}\|_2 \le C(q,\gamma)\left(\sqrt{\frac{s\|x^{\ast}\|_1\log p}{n}} + \sigma \sqrt{\frac{s\log p}{n}}\right),
\end{equation*}
where $C(q,\gamma)$ is a constant that depends only on $q$ and $\gamma$. Moreover, there exists a constant $C_1$ for which the following holds. In the case that $\sigma$ is unknown there is an estimator $\widehat{V}$ that satisfies with probability at least $1-n^{-2}$,
\begin{equation*}
    \max\left\{\left|\frac{\widehat{V}}{\sigma^2+\|x^{\ast}\|_1}-1\right|, \left|\frac{\sigma^2+\|x^{\ast}\|_1}{\widehat{V}}-1\right| \right\} \le C_1\sqrt{\frac{\log n}{\theta n}}.
\end{equation*}
In particular, setting $\widehat{\nu}:=\widehat{V}/\widehat{\mu}-1$ we have that 
\begin{equation*}
    \widehat{\nu} \rightarrow \frac{\sigma^2}{\|x^{\ast}\|_1},
\end{equation*}
almost surely provided that $\theta n \gg \log n$.
\end{proposition}

We add two remarks about our result. The first one is that the additional term $\sigma\sqrt{s\log p/n}$ in the convergence rate of the LASSO solution $\|\widehat{x}-x^{\ast}\|_2$, matches the standard convergence rate of the LASSO without Poisson noise \eqref{eq:convergence)_rate_standard_LASSO}. Therefore, it is sharp up to an absolute constant. The second one is that $\widehat{V}\sqrt{8\log p/n} = \sigma\sqrt{ 8\log p/n} +o(d_{old})$, thus the knowledge of $\sigma$ is not needed for the choice of $\lambda$.

Recall that using the design matrix $A_{new}$ we obtain $\theta n$ identically distributed samples $y_j$ given by \eqref{eq:observations_i.i.d_poisson_gauss} whose mean is $\|x^{\ast}\|_1$ and variance is $\|x^{\ast}\|_1 + \sigma^2$. Since one-dimensional variance estimation is also a mean estimation problem, the construction of $\widehat{V}$ is identical to $\widehat{\mu}$, so we leave the proof of the second part of Proposition \ref{main_prop_consistency} to the Appendix. Thanks to Proposition \ref{prop:estimation_l1norm_x*}, we can estimate $\|x^{\ast}\|_1$ accurately (if it is unknown). Next, one can immediately  estimate $\sigma^2/\|x^{\ast}\|_1$ as follows. Let $\widehat{\nu}$ as above and notice that almost surely
\begin{equation*}
\frac{\widehat{V}}{\widehat{\mu}} = \frac{\widehat{V}}{\sigma^2+\|x^{\ast}\|_1}\frac{\sigma^2+\|x^{\ast}\|_1}{\widehat{\mu}}  \rightarrow \frac{\sigma^2}{\|x^{\ast}\|_1} + 1,
\end{equation*}
therefore $\widehat{\nu}$ converges almost surely to $\sigma^2/\|x^{\ast}\|_1$. 

On the other hand, it might be of interest (at least theoretically) to derive guarantees for the scaled LASSO under Poisson noise as it does not require to modify the design matrix $A$. Motivated by this fact, we also show how to extend the scaled LASSO \eqref{eq:scaled_LASSO} to the Poisson-Gauss model.

The rest of this chapter is divided into two sections. The first one finishes the first part proof of Proposition \ref{main_prop_consistency} that establishes that $\|\widehat{x}-x^{\ast}\|_2/\sqrt{\|x^{\ast}\|_1}=o(1)$, i.e, the consistency of LASSO when $y$ is given by \eqref{eq:main_regression_model}. The second one analyzes the consistency of the scaled LASSO under the same model.

\subsection{Consistency of the LASSO}
Our starting point is to recall some machinery from \cite*{hunt2018data} to extend some of their results to the case of additional additive Gaussian noise. To this end, we first recall the definition of the well-known restricted eigenvalue condition that characterizes which design matrices are suitable for the LASSO algorithm.
\begin{definition}[Restricted Eigenvalue Condition]
\label{def:REC}
A matrix $\Tilde{A}$ satisfies the restricted eigenvalue condition $REC(\kappa_1.\kappa_2)$ if there exists $\kappa_1,\kappa_2>0$ such that for every $x\in \mathbb{R}^p$,
\begin{equation*}
    \|\widetilde{A}x\|_2 \ge \kappa_2 \|x\|_2 - \kappa_1 \|x\|_1.
\end{equation*}
\end{definition}
Also, an important term to quantify the performance of the LASSO is the following coherency-type parameter. Recall the choice of the regularization parameter in Theorem \ref{thm:consistency_LASSO_Poisson_only} was  $\lambda = \gamma d$, where $\gamma$ is an absolute constant. In what follows, we take $d$ to satisfy that 
\begin{equation}
\label{ineq:assumption_on_d}
    \|\widetilde{A}^T(\widetilde{y}-\widetilde{A}x^{\ast})\|_{\infty}\le d.
\end{equation}

The following result is central to our analysis. 
\begin{proposition}
\label{prop:performance_LASSO}
Fix $\varepsilon>0$. Suppose that $\Tilde{A}$ satisfies the restricted eigenvalue condition $REC(\kappa_1,\kappa_2)$ and that \eqref{ineq:assumption_on_d} holds for a certain $d>0$. For any $\gamma>2$, then there exists an absolute constant $c>0$ and a constant $\rho(\gamma)$ such that for any set $S\subset [p]$ satisfying that
\begin{equation}
\label{ineq:conditon_card_S}
    |S|\le \left(\frac{\kappa_2-\varepsilon}{\kappa_1\rho(\gamma)}\right)^2,
\end{equation}
the solution $\widehat{x}(\Tilde{A},\gamma d,y)$ of the LASSO program \ref{eq:LASSO_optimization} satisfies
\begin{equation*}
    \|\widehat{x}-x^{\ast}\|_2 \le c\left( \frac{\rho(\gamma)}{\varepsilon^2}\sqrt{s}d\right).
\end{equation*}
\end{proposition}
We are ready to prove the first part of Proposition \ref{main_prop_consistency}.
\begin{proof}

Notice that the REC condition and the inequality \eqref{ineq:conditon_card_S} do not depend on the noisy vector $g$, therefore it suffices to recall \cite*[Proposition 2]{hunt2018data} to ensure that there are absolute constants $c,c'>0$ such that for
\begin{equation}
\label{eq:kappas_REcondition}
\kappa_1 = \frac{c}{q(1-q)}\sqrt{\frac{\log p}{n}} \quad \kappa_2 = \frac{1}{4},
\end{equation}
the random matrix $\widetilde{A}$ \eqref{eq:tilde_A} satisfies the $REC(\kappa_1,\kappa_2)$ condition with probability at least $1-e^{-c'n}$. 

Next, we proceed to estimate $d$ in \eqref{ineq:assumption_on_d}, but we require some notation first. Recall that by definition
\begin{equation*}
\widetilde{y}:= \frac{1}{(n-1)\sqrt{nq(1-q)}} (ny - \sum_{l=1}^n y_l \mathbf{1}_{n\times 1}),
\end{equation*}
where $y_l = \mathcal{P}((Ax^{\ast})_l)+g_l$ and $y=(y_1,\ldots,y_n)$. To compare with the previous model, we adopt the convention $y_l^{old} = \mathcal{P}((Ax^{\ast})_l)$, $y^{old}=(y_1^{old},\ldots,y_n^{old})$ and $\widetilde{y}^{old}$ is defined analogously. The important feature of the Gaussian noise is that it is additive, which allows us to break the analysis as follows
\begin{equation*}
\widetilde{y} = \widetilde{y}^{old} + \underbrace{\frac{\sigma}{(n-1)\sqrt{nq(1-q)}} (ng - \sum_{l=1}^n g_l \mathbf{1}_{n\times 1})}_{:=G},
\end{equation*}
and by the triangle inequality
\begin{equation*}
 \|\widetilde{A}^T(\widetilde{y}-\widetilde{A}x^{\ast})\|_{\infty} \le  \| \widetilde{A}^T(\widetilde{y}^{old}-\widetilde{A}x^{\ast})\|_{\infty} + \|\widetilde{A}^TG\|_{\infty}.
\end{equation*}
The first term on the right-hand side was already studied in \cite*{hunt2018data}; in particular, we have that with high-probability
\begin{equation*}
\| \widetilde{A}^T(\widetilde{y}^{old}-\widetilde{A}x^{\ast})\|_{\infty} \lesssim \sqrt{\frac{\log p\|x^{\ast}\|_1}{nq}}  + \frac{\log p}{nq}.
\end{equation*}
Notice that the first term dominates when $nq>\log p$. All that remains is to deal with the noisy term that depends solely on $\widetilde{A}$ and $g$, but not on $\widetilde{y}$. Denote the first row of $\widetilde{A}$ by $\widetilde{a}_1^T$, and consider the first entry of $\widetilde{A}^T(G)$, namely
\begin{equation*}
\begin{split}
&|\widetilde{A}^T(G)|_1 = \frac{\sigma}{\sqrt{n}}\sum_{j=1}^n \frac{1}{\sqrt{q(1-q)}}(\widetilde{a}_1^T)_{j} (g_j +\frac{1}{n-1}\sum_{l\neq j}g_l)\\
&= \frac{\sigma}{\sqrt{n}}\sum_{j=1}^n \frac{1}{\sqrt{q(1-q)}}\widetilde{a}_{1j}g_j +  \frac{1}{\sqrt{n}}\sum_{j=1}^n \frac{\sigma}{\sqrt{q(1-q)}}\widetilde{a}_{1j}\frac{1}{n-1}\sum_{l\neq j}g_l\\
&=\frac{\sigma}{n\sqrt{q(1-q)}}\sum_{j=1}^n \frac{1}{\sqrt{q(1-q)}}(a_{1j}-q)g_j\\
&+  \frac{\sigma}{n\sqrt{q(1-q)}}\sum_{j=1}^n \frac{1}{\sqrt{q(1-q)}}(a_{1j}-q)\frac{1}{n-1}\sum_{l\neq j}g_l\\
&:= Z_1 + Z_2.
\end{split}
\end{equation*}
To estimate $Z_1$, notice that by independence between $A$ and $g$, we may condition on $A$ to obtain that
\begin{equation*}
\frac{1}{n}\sum_{j=1}^n \frac{\sigma}{\sqrt{q(1-q)}}(a_{1j}-q)g_j \sim N\left(0, \frac{\sigma^2}{n^2}\|\widetilde{a}_1^T\|_2^2\right).
\end{equation*}
Also, by Hoeffiding's inequality, there is an event $\Omega$ for which
\begin{equation*}
    \mathbb{P}\left(| \ \|\widetilde{a}_1^T\|_2^2-n \ |\ge \sqrt{\frac{n\log p}{q(1-q)}}\right)\le 2e^{-2\log p}=\frac{2}{p^2}.
\end{equation*}
On that event, using that $n\ge \log p$ it follows immediately that with probability at least $1-2e^{-t^2/2}$,
\begin{equation*}
\frac{1}{n}\sum_{j=1}^n \frac{\sigma}{\sqrt{q(1-q)}}(a_{1j}-q)g_j \le \sigma \sqrt{\frac{c_q}{n}}t,
\end{equation*}
where $c_q= (1+\sqrt{1/q(1-q)})$. It suffices to choose $t = 4\log p$ to obtain the probability that
\begin{equation*}
\frac{1}{n}\sum_{j=1}^n \frac{\sigma}{\sqrt{q(1-q)}}(a_{1j}-q)g_j \le \sigma \sqrt{\frac{4c_q\log p}{n}}
\end{equation*}
is at most $2p^{-2}$. Finally, notice that
\begin{equation*}
\begin{split}
\mathbb{P}\left(|Z_1|\ge \sigma \sqrt{\frac{4c_q\log p}{q(1-q)n}} \right) &\le \mathbb{P}\left( |Z_1| \ge \sigma \sqrt{\frac{4c_q\log p}{q(1-q)n}} |\Omega\right)\mathbb{P}(\Omega) + \mathbb{P}(\Omega^c) \\
&\le \frac{2}{p^2}+ \frac{2}{p^2}=\frac{4}{p^2},
\end{split}
\end{equation*}
therefore on that event
\begin{equation*}
\frac{1}{n\sqrt{q(1-q)}}\sum_{j=1}^n \frac{\sigma}{\sqrt{q(1-q)}}(a_{1j}-q)g_j \le \sigma \sqrt{\frac{4c_q\log p}{q(1-q)n}}:=\sqrt{\frac{4c_q'\log p}{n}}.
\end{equation*}
To estimate $Z_2$, notice that for every $j\in [n]$,
\begin{equation*}
    g(j)':=\frac{\sigma}{(n-1)}\sum_{l\neq j} g_l \sim N\left(0,\frac{\sigma^2}{n-1}\right).
\end{equation*}
Thus, we have that $g(j)' \sim \sigma/\sqrt{n-1} g$. Recall the standard Gaussian tail inequality
\begin{equation*}
\mathbb{P}(g\ge t) \le 2 e^{-t^2/2},
\end{equation*}
and choose $t=\sqrt{8\log p}$ to obtain that with probability at least $1-2/p^4$,
\begin{equation*}
    |g(j)'|\le \frac{\sigma}{\sqrt{n-1}}\sqrt{8\log p}.
\end{equation*}
By union bound over all $j\in [n]$, we ensure that with probability at least $1-n/p^4$, for every $j\in [n]$
\begin{equation*}
    \left|\frac{\sigma}{(n-1)}\sum_{l\neq j}g_l\right| \le \sigma \sqrt{\frac{8\log p}{n}}.
\end{equation*}
On that event (which we denote by $\Omega_1$), we have that
\begin{equation*}
    |Z_2| \le \sigma \sqrt{\frac{8\log p}{n}} \frac{1}{\sqrt{q(1-q)}} \frac{1}{n}\sum_{j=1}^n \frac{\sigma}{\sqrt{q(1-q)}}(a_{1j}-q).
\end{equation*}
By independence between $A$ and $g$, once we condition on $g$, we still have that $a_{1j}-q/\sqrt{q(1-q)}$ are independent mean zero random variables with variance one. Thus, by Bernstein's inequality, with probability at least $1-e^{-2nq(1-q)}$
\begin{equation*}
\frac{1}{\sqrt{q(1-q)}}\left|\frac{1}{n}\sum_{j=1}^n \frac{1}{\sqrt{q(1-q)}}(a_{1j}-q)\right| \le \frac{1}{\sqrt{q(1-q)}}2\sqrt{q(1-q)} = 2.
\end{equation*}
Finally, 
\begin{equation*}
\begin{split}
\mathbb{P}\left(|Z_2|\ge \sigma \sqrt{\frac{8\log p}{n}} \right) &\le \mathbb{P}\left(|Z_2|\ge \sigma \sqrt{\frac{8\log p}{n}}|\Omega_1 \right) \mathbb{P}(\Omega_1) + \mathbb{P}(\Omega_1^c) \\
&\le e^{-nq(1-q)/2}+\frac{2n}{p^4}\le e^{-2nq(1-q)} + \frac{2}{p^3}.
\end{split}
\end{equation*}
The same estimate $|\widetilde{A}^T(G)|_k =O(\sigma\sqrt{\log p/n})$ holds simultaneously for every $k\in [p]$ by union bound. We conclude that on a high-probability event, it is enough to choose $d$ of order
\begin{equation*}
d \simeq \sqrt{\frac{\log p\|x^{\ast}\|_1}{nq}} + \frac{\log p \|x^{\ast}\|_1}{n} + \frac{\log p}{nq} + \sigma\sqrt{\frac{c'(q)\log p}{n}}.
\end{equation*}
Recall that we assume that $\|x^{\ast}\|_1$ and $\sigma$ are known, then $d$ is also known. Applying Proposition \ref{prop:performance_LASSO} we obtain that with high-probability
\begin{equation*}
\|\widehat{x}-x^{\ast}\|_2 \lesssim \sqrt{\frac{s\|x^{\ast}\|_1\log p}{nq}} + \sigma\sqrt{\frac{c'(q)s\log p}{n}}.
\end{equation*}
\end{proof}
\subsection{Consistency of the Scaled LASSO}
In this section, we focus on the estimation of the standard deviation of the additive Gaussian noise, namely $\sigma$ in \eqref{eq:main_regression_model} using the scaled LASSO. Readers who are only interested in the proof of our main result that is Theorem \ref{thm:main_result} may skip this section. To start, we recall that the scaled LASSO estimator is given by \eqref{eq:scaled_LASSO}. The main result is the following.
\begin{proposition}
\label{prop:consistency_scaled_Lasso}
Consider the same setup as in Proposition \ref{main_prop_consistency}. Then the solution $\widehat{\sigma}(\lambda')$ of \eqref{eq:scaled_LASSO},  for the choice of 
\begin{equation*}
\lambda'=\sqrt{\frac{\|x^{\ast}\|_1\log p}{n}},
\end{equation*}
satisfies with high probability
\begin{equation*}
    \max\left\{\left|\frac{\widehat{\sigma}}{\sigma}-1\right|, \left|\frac{\sigma}{\widehat{\sigma}}-1\right| \right\} \le K\sqrt{\|x^{\ast}\|_1\frac{s\log p}{n}},
\end{equation*}
where $K$ is an absolute constant.
\end{proposition}

The rest of this section is dedicated to the proof of the Proposition \ref{prop:consistency_scaled_Lasso}. To start, we recall some background from the theory of the LASSO algorithm. We refer the reader to \cite*{buhlmann2011statistics} for a more detailed presentation. Let $\mathcal{C}(\xi, T)$ be the cone 
\begin{equation*}
\mathcal{C}(\xi, T):=\{x\in \mathbb{R}^p: \|x_{T^c}\|_1\le \xi \|x_{T}\|_1 \},
\end{equation*}
and define the compatibility factor as
\begin{equation*}
\kappa(\xi,s):= \inf_{x\in \mathcal{C}(\xi, T^{\ast})/\{0\}} \sqrt{s}\frac{\|\widetilde{A}x\|_2}{\|x_{T^{\ast}}\|_1},
\end{equation*}
where $T^{\ast}$ is the support of the best $s$-term approximation of $x$. The oracle prediction error bound is defined by 
\begin{equation}
\label{eq:oracle_prediction_error}
\eta(\xi,\lambda):= \sqrt{s}\frac{2\xi\lambda }{(\xi+1)\kappa(\xi,s)},
\end{equation}
and the oracle noise by 
\begin{equation}
\label{eq:oracle_noise}
    \sigma^{\ast}:= \frac{\|\widetilde{y}-\widetilde{A}x^{\ast}\|_2}{\sqrt{n}}.
\end{equation}
We are now ready to state the main auxiliary result due to \cite*[Theorem 1]{sun2012scaled}.
\begin{theorem}
Consider $\eta=\eta(\xi,\lambda)$ as in \eqref{eq:oracle_prediction_error}. Suppose there is a $\xi>1$ for which the following holds
\begin{equation}
\label{condition_scaled_lasso}
    \|\widetilde{A}^T(\widetilde{y}-\widetilde{A}x^{\ast})\|_{\infty} \le \sigma^{\ast}(1-\eta)\lambda\frac{(\xi-1)}{(\xi+1)}.
\end{equation}
Then \eqref{eq:scaled_LASSO} outputs $\widehat{\sigma}(\lambda)$ that satisfies
\begin{equation*}
    \max\left\{\left|\frac{\widehat{\sigma}}{\sigma}-1\right|, \left|\frac{\sigma}{\widehat{\sigma}}-1\right| \right\} \le \eta.
\end{equation*}

\end{theorem}
Also, we need the following technical result known as Anderson's lemma \cite*{le2000asymptotics}. 
\begin{lemma}
\label{lemm:Anderson}
Let $G\sim N(0,\Sigma)$ be a multivariate Gaussian. Then for any norm $\|\cdot\|$ on $\mathbb{R}^p$, the following holds
\begin{equation*}
    \inf_{x\in \mathbb{R}^d}\mathbb{E}\|G+x\|= \mathbb{E}\|G\|
\end{equation*}
\end{lemma}
We are now ready to prove Proposition \ref{prop:consistency_scaled_Lasso}.
\begin{proof}
Recall that $\widetilde{A}$ satisfies the $REC(\kappa_1,\kappa_2)$ \eqref{eq:kappas_REcondition} (with high-probability). Therefore, it immediately follows that
\begin{equation*}
    \kappa(\xi,s) \ge \frac{1}{4} - \frac{(1+\xi)c}{q(1-q)}\sqrt{\frac{s\log p}{n}}.
\end{equation*}
Thus, assuming that $n\ge cs\log p$ for some well-chosen $c>0$, we ensure that $\kappa(\xi,s)\ge c_1>0$. It is straightforward to verify that for
\begin{equation*}
    \lambda = \sqrt{\frac{\log p\|x^{\ast}\|_1}{n}},
\end{equation*}
we obtain that
\begin{equation*}
    \eta \le 2c_1 \frac{\xi}{(1+\xi)}\sqrt{\|x^{\ast}\|_1\frac{s\log p}{n}}.
\end{equation*}
Next, recall that from the previous section, we have that on a high-probability event,
\begin{equation*}
\|\widetilde{A}^T(\widetilde{y}-\widetilde{A}x^{\ast})\|_{\infty} \lesssim \sqrt{\frac{\log p\|x^{\ast}\|_1}{nq}}  + \sigma\sqrt{\frac{c'(q)\log p}{n}}.
\end{equation*}
Recall that $\sigma^{\ast}$ is the oracle noise given by \eqref{eq:oracle_noise}. Therefore to verify that \eqref{condition_scaled_lasso} holds, it suffices to check that
\begin{equation*}
    \sigma \lesssim \sigma^{\ast}.
\end{equation*}
This implies that $\widehat{\sigma}$ converges to $\sigma$ as required. The rest of this section is dedicated to showing that $\sigma^{\ast}\ge \sigma/4$ with high probability. In fact, recall that
\begin{equation*}
G=
\begin{pmatrix}
g_1 + \frac{1}{(n-1)}\sum_{i\neq 1}g_i\\
\vdots\\
g_n + \frac{1}{(n-1)}\sum_{i\neq n}g_i
\end{pmatrix},
\end{equation*}
and notice that by independent between $y_1,\ldots,y_n$ and $g$, we have that
\begin{equation*}
\sigma^{\ast} = \frac{\|\widetilde{y}^{old}-\widetilde{A}x^{\ast} + G\|_2}{\sqrt{n}}:=\frac{\|G+z\|_2}{\sqrt{n}},
\end{equation*}
for some $z$ that is independent from $G$. We first show that $G$ is a joint Gaussian vector. To verify this, it suffices to show that for every $v \in S^{n-1}$, $\langle G,v\rangle$ is Gaussian (see \cite*[Exercise 3.3.4]{vershynin2018high} for a formal proof of this standard fact). Since $G$ has independent Gaussian entries, we just need to check that the following random variable
\begin{equation*}
    \sum_{i=1}^n G_iv_i = \sigma\sum_{i=1}^n g_i ( v_i + \frac{1}{n-1}\sum_{l\neq i}v_l),
\end{equation*} 
is indeed Gaussian. Since it is a linear combination of independent Gaussians, such a random variable fails to be Gaussian if and only if for every $j\in [n]$
\begin{equation*}
    v_i + \frac{1}{n-1}\sum_{l\neq i}v_l =0 \quad \text{equivalently} \quad (n-2)v_i = -\sum_{l=1}^n v_l.
\end{equation*}
It follows that this only holds when $(n-2)v_1=\ldots=(n-2)v_n=0$. Thus, either $n=2$ (we rule out this case by assuming that $n\ge 3$) or $v_1=\ldots=v_n=0$, but the latter is impossible because $v \in S^{n-1}$. The fact that $G$ is a Gaussian vector follows.

Next, since $G$ is a Gaussian vector and it is independent of $\widetilde{y}^{old}$, we have that conditioned on $\widetilde{y}$, the vector $G+z$ still distributed as a Gaussian with mean $z$. Setting  $\Sigma=\mathbb{E}GG^{T}$ to be the covariance matrix of $G$ and applying the Lipschitz concentration inequality for Gaussians  \cite*[Theorem 5.2.2]{vershynin2018high} with the choice of a Lipschitz function  $f(x)=\Sigma^{1/2}x+z$, it follows that with probability at least $1-e^{-nt^2}$
\begin{equation}
\label{ineq:lipschitz_Gaussian_part1}
\frac{1}{\sqrt{n}}\|G+z\|_2 \ge \mathbb{E} \frac{1}{\sqrt{n}}\|G+z\|_2 - t\sqrt{\|\Sigma\|}.
\end{equation}
Next, we bound the quantity $\|\Sigma\| $. It follows that
\begin{equation*}
    \Sigma_{ii} = \sigma^2(1+\frac{1}{n-1}) \quad \text{and for $i\neq j$} \quad \Sigma_{i,j}=2\frac{\sigma^2}{(n-1)} + \sigma^2\frac{n-2}{(n-1)^2}.
\end{equation*}
Thus, by Gershgorin's circle theorem
\begin{equation*}
    \left|\|\Sigma\| -\sigma^2(1+\frac{1}{n-1}) \right| \le 3\sigma^2,
\end{equation*}
and then $\|\Sigma\|\le 5\sigma^2$. From \eqref{ineq:lipschitz_Gaussian_part1}, we obtain that with probability $1-e^{-n/40}$
\begin{equation}
\label{ineq:lipschitz_Gaussian_part2}
\frac{1}{\sqrt{n}}\|G+z\|_2 \ge \mathbb{E} \frac{1}{\sqrt{n}}\|G+z\|_2 - \frac{\sigma}{2}.
\end{equation}
It suffices to show that $\mathbb{E}\|G+z\|_2 \ge (1-o(1))\sigma \sqrt{n}$. More accurately, we apply Anderson's lemma \ref{lemm:Anderson} to obtain that $$\mathbb{E}\|G+z\|_2 \ge \mathbb{E}\|G\|_2.$$ Next, it is straightforward to verify that
\begin{equation*}
    \mathbb{E}\|G\|_2^2 \ge \sigma^2 n.
\end{equation*}
Secondly, by applying again the Lipschitz concentration inequality
\begin{equation*}
\Var(\|G\|_2) = \int_{0}^{\infty}\mathbb{P}(\left|\|G\|_2-\mathbb{E}\|G\|_2\right|^2 \ge t) dt\le 2\sigma^2\int_{0}^{\infty}e^{-t/2} = 4\sigma^2.
\end{equation*}
using the fact that
\begin{equation*}
  \mathbb{E}\|G\|_2 \ge \sqrt{\mathbb{E}\|G\|_2^2 -  \Var(\|G\|_2)},
\end{equation*}
we obtain that
\begin{equation}
\label{ineq:expectation_Gaussian}
    \frac{1}{\sqrt{n}}\mathbb{E}\|G+z\|_2 \ge \frac{1}{\sqrt{n}}\mathbb{E}\|G\|_2 \ge \sigma\left(1-\frac{2}{\sqrt{n}}\right).
\end{equation}
Finally, by the law of iterated expectation, we take expectation with respect to $z$ and apply \eqref{ineq:lipschitz_Gaussian_part2} and \eqref{ineq:expectation_Gaussian} to obtain that with high probability 
\begin{equation*}
    \sigma^{\ast} \ge \frac{1}{4}\sigma.
\end{equation*}
The constant $1/4$ follows by assuming that $n\ge 16$.
\end{proof}

\section{The Debiased Algorithm under the Poisson-Gauss noise}
\label{sec:main_result}
In this section, we put together the arguments from Section \ref{sec:UQ_poisson_alone} and \ref{sec:consistency_poissongauss} to provide the final form of our estimator. Without loss of generality, assume that $(1-\theta)n$ is an integer. Recall that the bias $B$ is given by \eqref{eq:expression_known_bias} and that the modified design matrix $A_{new}$ is given by \eqref{eq:new_design_matrix}.

Notice that, in general, $B$ is unknown  because it depends on $\|x^{\ast}\|_1$, however we bypass this obstacle by estimating $\|x^{\ast}\|_1$ accurately using Proposition \ref{prop:estimation_l1norm_x*}. The final algorithm is the following.
\begin{algorithm}[H]
\caption{Debiased Algorithm under Poisson-Gauss noise}
\begin{algorithmic}
\Require Design matrix $ \begin{pmatrix}
A\\
A_{new}'
\end{pmatrix} \in \{0,1\}^{2n\times p}$, a response vector $(y,y')\in \mathbb{R}^{2n}$ from \eqref{eq:main_regression_model}. 
\State
\State Compute $\widehat{\mu}$ using the last $\theta n$ rows of $A_{new}'$ via Proposition \ref{prop:estimation_l1norm_x*}.

\State Compute $\Tilde{A}$, $\Tilde{A'}_{new}$ by \eqref{eq:tilde_A} and  $\Tilde{y}$,$\Tilde{y}'$ by \eqref{eq:tilde_y}.

\State Estimate $\sigma$ using $\Tilde{A}$ and $\Tilde{y}$ via Proposition \ref{main_prop_consistency} (if necessary).

\State Compute the LASSO solution $\widehat{x}$ as in Proposition \ref{main_prop_consistency} using the first $(1-\theta)n$ rows of $\Tilde{A'}_{new}$ and $\Tilde{y}'$.

\State Compute the estimated bias $\widehat{B}$ replacing $\|x^{\ast}\|_1$ by $\widehat{\mu}$ in \eqref{eq:expression_known_bias}.
\\

\Return
\begin{equation}
    \widehat{x}^{d}:= \widehat{x} + \frac{1}{\sqrt{nq(1-q)}}\Tilde{A}^T(y-A\widehat{x}) - \frac{1}{\sqrt{n}}\widehat{B}.
\end{equation}
\end{algorithmic}
\end{algorithm}

The proof of Theorem \ref{thm:main_result} readily follows from the estimator above and the arguments presented in the text. Notice that we may replace $n$ by $(1-\theta)n$ in the estimates for $\eta,\Delta,\eta'$ \eqref{eq:decomposition_deabised_final}. We remark that it does not affect the conclusion of the statement, not even the constants, because we are allowed to take a vanishing $\theta$ if necessary.

\section{Numerical Experiments}
\label{sec:numerical}
In this section, we perform some numerical experiments with synthetic data to support our theory. In this context, we generate the data according to \eqref{eq:main_regression_model} with $p=2\cdot 10^4$, $s=100$ and $n=\lceil s \log^2p \rceil$. The ground truth $x^{\ast}$ is generated at random. We perform the experiments using the MATLAB package of LASSO.

For this experiment, we computed the confidence levels via $\eqref{eq:confidencen_levels}$ by setting $\alpha=0.1$. Mathematically, it means that the probability for a sample to be outside of the estimated interval is asymptotically equal to $0.1$ (see \eqref{eq:asymptotically_valid_CI}). From a practical perspective, data splitting is not convenient, so we run our main algorithm ten times without splitting the data (it only makes the performance worse). As we shall see, data splitting does not seem necessary for the de-biased LASSO.

Moreover, we first assumed the knowledge of $\|x^{\ast}\|_1$ and $\sigma^2$. For the regularization parameter, we used $d_{old}$ given by \eqref{eq:easy_choice_regularization}. We tested different values $\gamma$ ranging from $2.01$ to $70$ (when it starts to output meaningless results), each of them seems to perform equally well on the $s$ non-zero coefficients, however for $\gamma=50$ the LASSO seems to achieve the best error on the zero coefficients. We avoided using cross-validation as we did not analyze it theoretically, and it is computationally expensive. We ran our algorithm ten times, and for each run, we computed the number of mistakes, namely the number of samples that were outside of the estimated confidence interval.

The average number (over ten trials) of mistakes was $3$ and the largest number among those was $6$ mistakes. These experimental numbers are in agreement with the theory as roughly one should expect at most $10$ mistakes. Indeed, for each sample, the probability of being outside of the estimated interval is $0.1$, and we have $s=100$ non-zero coefficients. It seems that our estimates on the number of samples are too conservative.

The figure below illustrates the coverage corresponding to the worst-case obtained among ten trials. The coefficients are re-arranged in a non-increasing order to facilitate the visualization. Clearly, this does not affect the number of mistakes.
\begin{figure}
\centering
\includegraphics[width=0.6\textwidth]{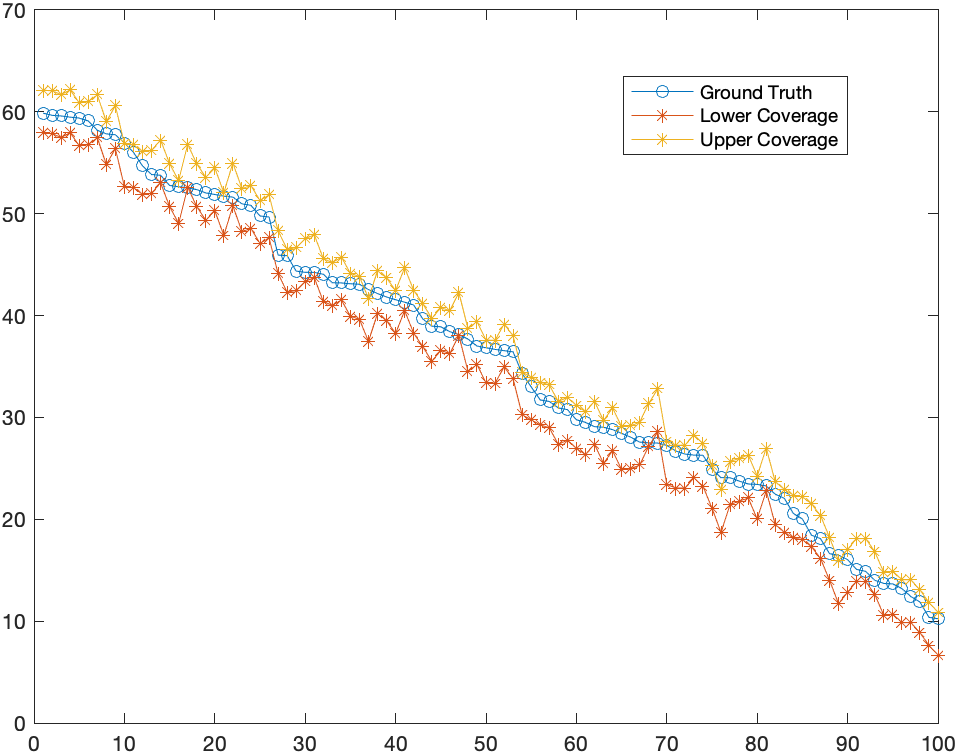}
\medskip
\caption{Confidence intervals for the $100$-largest coefficients arranged in the non-increasing order.}
\end{figure}
As one can easily see, even the samples that are outside of the estimate interval are quite close to its interval, meaning that the errors on those are negligible.

Next, we run some experiments on the estimation of $\|x^{\ast}\|_1$ and $\nu:=\sigma^2/\|x^{\ast}\|_1$. Due to the light tail of Poisson and Gaussian random variables, the empirical mean already provides accurate estimates. Indeed, we only used $100$ samples among the total of $9808$ samples to estimate those parameters. To evaluated the error on estimating $\|x^{\ast}\|_1$, we computed the following error
\begin{equation*}
   \varepsilon_{mean}:=\max\left\{\left|\frac{\|x^{\ast}\|_1}{\widehat{\mu}}-1\right|,\left|\frac{\widehat{\mu}}{\|x^{\ast}\|_1}-1\right|\right\} ,
\end{equation*}
and to estimate $\nu$, we used  
\begin{equation*}
   \varepsilon_{var}:=\max\left\{\left|\frac{\nu}{\widehat{\nu}}-1\right|,\left|\frac{\widehat{\nu}}{\nu}-1\right|\right\}.
\end{equation*}
After running ten times, we obtained that the largest value of $\varepsilon_{mean}$ and $\varepsilon_{var}$ was $4.5\cdot 10^{-3}$ and $0.19$, respectively. The average value of $\varepsilon_{mean}$ and $\varepsilon_{var}$ was $1.4 \cdot 10^{-3}$ and $0.17$, respectively. Due to such small errors, the uncertainty with respect to $\|x^{\ast}\|_1$ and $\nu$ is negligible.


\newpage
\section*{Appendix}
The goal of this Appendix is to prove Proposition \ref{prop:estimation_l1norm_x*} and the second part of Proposition \ref{main_prop_consistency}. To this end, we invoke the standard median of means procedure to construct an estimator $\widehat{\mu}$ \cite*[Theorem 5.2]{lugosi2017lectures}.
\begin{proposition}
\label{prop:mom_1d}
Let $X_1,\ldots, X_n$ be $n$ independent copies of a random variable $X$ with mean $\mu$ and variance $\sigma^2$. Then for any $t>0$, there exists an estimator $\widehat{\mu}$ satisfying with probability at least $1-e^{-t}$,
\begin{equation*}
    |\widehat{\mu}-\mu|\le \sigma \sqrt{\frac{32 t}{N}}.
\end{equation*}
\end{proposition}
The estimator proceeds as follows: it splits the sample into $k=\lfloor8t \rfloor$ blocks $B_1,\ldots, B_{k}$ of equal size each and computes 
\begin{equation*}
    \widehat{\mu}_j  = \frac{1}{|B_j|}\sum_{i\in B_j} X_i,
\end{equation*}
then it outputs $\widehat{\mu}$ defined to be the median of $\widehat{\mu}_1,\ldots,\widehat{\mu}_k$. Clearly, it can be computed in polynomial time.

Next, recall from \eqref{eq:observations_i.i.d_poisson_gauss} that the observations that we collect by modifying the design matrix using $A_{new}$ are $\theta n$ i.i.d samples given by
\begin{equation*}
    y_j = \mathcal{P}_{j}(\|x^{\ast}\|_1) + \sigma g_j,
\end{equation*}
whose mean is $\|x^{\ast}\|_1$ and the variance is $\|x^{\ast}\|_1+\sigma^2$. It follows immediately from Proposition \ref{prop:mom_1d} that there is an estimator $\widehat{\mu}$ that satisfies with probability at least $1-e^{-t}$,
\begin{equation*}
|\widehat{\mu}-\|x^{\ast}\|_1| \le (\sqrt{\|x^{\ast}\|_1}+\sigma)\sqrt{\frac{32t}{\theta n}},
\end{equation*}
thus setting $t=2\log n$, we obtain that with probability at least $1-n^{-2}$
\begin{equation*}
|\widehat{\mu}-\|x^{\ast}\|_1| \le 8(\sqrt{\|x^{\ast}\|_1}+\sigma)\sqrt{\frac{\log n}{\theta n}}.
\end{equation*}

The proof of the second part of Proposition \ref{main_prop_consistency} follows the same proof. Indeed, by symmetrizing the samples, precisely $\Delta y_1:=(y_1-y_2)/\sqrt{2},\ldots, \Delta y_{\theta n/2}:=(y_{\theta n -1}-y_{\theta n})/\sqrt{2}$ we may assume that we collect $\theta n/2$ independent samples distributed as $\Delta y$  whose mean zero and variance $\sigma^2$. Clearly, both Poisson and Gaussian distributions have finite fourth moments, and then another application of Proposition \ref{prop:mom_1d} to $(\Delta y)^2$ suffices to construct $\widehat{V}$.

\end{document}

%% file: arxiv_style.tex
\usepackage[letterpaper, left=1in, right=1in, top=1in, bottom=1in]{geometry}

\usepackage[utf8]{inputenc} 
\usepackage[T1]{fontenc}    
\usepackage{url}            
\usepackage{booktabs}       
\usepackage{amsfonts}       
\usepackage{nicefrac}       
\usepackage{microtype}      
\usepackage{authblk}
\usepackage{tocloft}            
\usepackage[nottoc]{tocbibind}
\usepackage{enumitem}
\usepackage{tcolorbox}
\usepackage{algorithm}
\usepackage{algpseudocode}
\usepackage{graphicx}
\usepackage{pstricks}
\usepackage{appendix}
\usepackage{dsfont}
\usepackage{amsmath}
\usepackage{amssymb}
\usepackage{amsthm}
\usepackage{comment}
\usepackage{tcolorbox}
\usepackage[nottoc]{tocbibind}

\usepackage{breakcites}

\usepackage{etoolbox}
\usepackage{comment}
\newtoggle{draft}
\togglefalse{draft}

\usepackage{color-edits}

\usepackage{mathrsfs}

\usepackage{algorithm}
\usepackage{verbatim}

\usepackage{multicol}

\usepackage{colortbl}
\usepackage{bbm}
\usepackage{setspace}

\usepackage{transparent}

\usepackage{inconsolata}
\usepackage[scaled=.90]{helvet}
\usepackage{xspace}

\usepackage{pifont}

\usepackage{tikz-cd}

    \usepackage{algorithm}
    \usepackage{authblk}
     \usepackage{algpseudocode}
    \usepackage{fullpage}
    \usepackage{natbib}
    \usepackage{amsthm}
    \usepackage{color}
    \usepackage{amsmath}
    \usepackage{amsfonts}
    \usepackage{xcolor}
    \colorlet{linkequation}{blue}
    \usepackage{framed}
    \usepackage{times}
    \usepackage{amsmath}
    \usepackage{amsthm}
    \usepackage{amssymb}
    \usepackage{wasysym}
    \usepackage{mathrsfs}
    \usepackage{bbm}
    \usepackage{geometry} 
    \usepackage{mathtools}

    \newcounter{rcnt}[section]

    \usepackage{tikz}
\usepackage{lipsum}

\newlength{\RoundedBoxWidth}
\newsavebox{\GrayRoundedBox}
   {\setlength{\RoundedBoxWidth}{\dimexpr#1}
    \begin{lrbox}{\GrayRoundedBox}
       \begin{minipage}{\RoundedBoxWidth}}%
   {   \end{minipage}
    \end{lrbox}
    \begin{center}
    \begin{tikzpicture}%
       \draw node[draw=black,fill=black!10,rounded corners,%
             inner sep=2ex,text width=\RoundedBoxWidth]%
             {\usebox{\GrayRoundedBox}};
    \end{tikzpicture}
    \end{center}}
    
    \def\argmin{\mathop{\rm argmin}}

    \newcommand{\Var}{\operatorname{Var}}

    \usepackage{lipsum}
    \usepackage{appendix}
    
    \usepackage[T1]{fontenc}
    \usepackage{authblk}

\PassOptionsToPackage{hypertexnames=false}{hyperref}  
\usepackage{parskip}

\usepackage[colorlinks=true, linkcolor=blue!70!black, citecolor=blue!70!black,urlcolor=black,breaklinks=true]{hyperref}
\usepackage{microtype}
\usepackage{hhline}

\makeatletter
\newcommand{\neutralize}[1]{\expandafter\let\csname c@#1\endcsname\count@}
\makeatother

\usepackage{algorithm}

\usepackage{natbib}
\bibpunct{(}{)}{;}{a}{,}{,}

\usepackage{amsthm}
\usepackage{mathtools}
\usepackage{amsmath}
\usepackage{bbm}
\usepackage{amsfonts}
\usepackage{amssymb}

\usepackage{xpatch}

\newtheorem{theorem}{Theorem}
\newtheorem*{theorem*}{Theorem}
\newtheorem{lemma}{Lemma}
\newtheorem{proposition}[theorem]{Proposition}

\newtheorem{definition}{Definition}
\theoremstyle{definition}
\newtheorem{example}{Example}

\theoremstyle{remark}

\AtBeginEnvironment{remark}{%
  \pushQED{\qed}%
}
\AtEndEnvironment{remark}{\popQED\endexample}

\makeatletter
  \renewenvironment{proof}[1][Proof]%
  {%
   \par\noindent{\bfseries\upshape {#1.}\ }%
  }%
  {\qed\newline}
  \makeatother

\xpatchcmd{\proof}{\itshape}{\normalfont\proofnameformat}{}{}
\newcommand{\proofnameformat}{\bfseries}

\usepackage[nameinlink,capitalize]{cleveref}

\crefformat{equation}{#2(#1)#3}
\Crefformat{equation}{#2(#1)#3}

\Crefformat{figure}{#2Figure #1#3}
\Crefname{assumption}{Assumption}{Assumptions}
\Crefformat{assumption}{#2Assumption #1#3}

\usepackage{crossreftools}
\usepackage{xparse}

\ExplSyntaxOn
\DeclareDocumentCommand{\XDeclarePairedDelimiter}{mm}
 {
  \__egreg_delimiter_clear_keys: 
  \keys_set:nn { egreg/delimiters } { #2 }
  \use:x 
   {
    \exp_not:n {\NewDocumentCommand{#1}{sO{}m} }
     {
      \exp_not:n { \IfBooleanTF{##1} }
       {
        \exp_not:N \egreg_paired_delimiter_expand:nnnn
         { \exp_not:V \l_egreg_delimiter_left_tl }
         { \exp_not:V \l_egreg_delimiter_right_tl }
         { \exp_not:n { ##3 } }
         { \exp_not:V \l_egreg_delimiter_subscript_tl }
       }
       {
        \exp_not:N \egreg_paired_delimiter_fixed:nnnnn 
         { \exp_not:n { ##2 } }
         { \exp_not:V \l_egreg_delimiter_left_tl }
         { \exp_not:V \l_egreg_delimiter_right_tl }
         { \exp_not:n { ##3 } }
         { \exp_not:V \l_egreg_delimiter_subscript_tl }
       }
     }
   }
 }

\keys_define:nn { egreg/delimiters }
 {
  left      .tl_set:N = \l_egreg_delimiter_left_tl,
  right     .tl_set:N = \l_egreg_delimiter_right_tl,
  subscript .tl_set:N = \l_egreg_delimiter_subscript_tl,
 }

\cs_new_protected:Npn \__egreg_delimiter_clear_keys:
 {
  \keys_set:nn { egreg/delimiters } { left=.,right=.,subscript={} }
 }

\cs_new_protected:Npn \egreg_paired_delimiter_expand:nnnn #1 #2 #3 #4
 {
  \mathopen{}
  \mathclose\c_group_begin_token
   \left#1
   #3
   \group_insert_after:N \c_group_end_token
   \right#2
   \tl_if_empty:nF {#4} { \c_math_subscript_token {#4} }
 }
\cs_new_protected:Npn \egreg_paired_delimiter_fixed:nnnnn #1 #2 #3 #4 #5
 {
  \mathopen{#1#2}#4\mathclose{#1#3}
  \tl_if_empty:nF {#5} { \c_math_subscript_token {#5} }
 }
\ExplSyntaxOff

\XDeclarePairedDelimiter{\supnorm}{
  left=\lVert,
  right=\rVert,
  subscript=\infty
  }